\documentclass[10pt, fleqn]{article}
\usepackage{graphicx}
\usepackage{latexsym}
\usepackage{xcolor}
\usepackage{amsmath,amsthm,amssymb}
\usepackage{mathrsfs}
\usepackage[hypcap]{caption}
\usepackage[all]{xy}
\usepackage[textwidth=14.5cm]{geometry}

\usepackage{tabularx}
\makeatletter
 
 \@addtoreset{equation}{section}
\makeatother
\makeatletter
\newcommand{\subjclass}[2][1991]{%
  \let\@oldtitle\@title%
  \gdef\@title{\@oldtitle\footnotetext{#1 \emph{Mathematics subject classification.} #2}}%
}
\newcommand{\keywords}[1]{%
  \let\@@oldtitle\@title%
  \gdef\@title{\@@oldtitle\footnotetext{\emph{Key words and phrases.} #1.}}%
}
\makeatother

\newcommand{\pf}{\textit{\sc Proof.} \ }
\newcommand{\remark}{\textit{Remark.} \ }

\theoremstyle{plain}
\newtheorem{df}{\sc \bf Definition}[section]
\newtheorem{thm}[df]{\sc \bf Theorem}

\newtheorem{prop}[df]{\sc \bf Proposition}
\newtheorem{lem}[df]{\sc \bf Lemma}
\newtheorem{cor}[df]{\sc \bf Corollary}

\theoremstyle{remark}

\def\X{{\mathfrak{X}(M)}} 
\def\dfrac#1#2{{\displaystyle\frac{#1}{#2}}}
\newcommand{\pari}{\frac{\partial}{\partial x_i}}
\newcommand{\parj}{\frac{\partial}{\partial x_j}}

\newcommand{\R}{\boldsymbol{R}}

\newcommand{\tr}{\mathrm{tr}}

\def\l{{\mathfrak{l}}}
\newcommand{\h}{\mathfrak{h}}

\newcommand{\g}{\mathfrak{g}}
\newcommand{\kk}{\mathfrak{k}}
\newcommand{\da}{\mathfrak{a}}

\newcommand{\sll}{\mathfrak{sl}}
\newcommand{\gl}{\mathfrak{gl}}

\newcommand{\cd}{\cdot}
 
\usepackage{enumerate}
\usepackage{multirow,bigdelim}
\usepackage{delarray}

\begin{document}
\title{Low dimensional Lie groups admitting left invariant flat projective or affine structures}
\author{Hironao Kato}
\date{\vspace{-2ex}}

\maketitle
\begin{abstract}
We prove that any real Lie group of dimension $\leq 5$ admits a left invariant flat projective structure. We also prove that a real Lie group $L$ of dimension $\leq 5$ admits a left invariant flat affine structure if and only if the Lie algebra of $L$ is not perfect. 
\end{abstract}
\section{Introduction} 
A left invariant flat affine structure (IFAS for short) is a 
torsion-free affine connection $\nabla$ on a Lie group $L$ such that $\nabla$ is left invariant and flat. A left invariant flat projective 
structure $[\nabla]$ (IFPS for short) is a projective 
equivalence class of an affine connection $\nabla$ on $L$ such that the left action 
of $L$ is a projective transformation and $\nabla$ is locally projectively equivalent to some flat affine connection (cf. Definition \ref{def 
of IFPS}). 

In this paper we consider the problem whether low dimensional Lie groups admit these geometric structures or not. Our main results are the following.
\begin{thm}\label{IFPS}
Any real Lie group of dimension $\leq 5$ admits an IFPS. 
\end{thm}

\begin{thm}\label{IFAC}
Let $L$ be a real Lie group of dimension $\leq 5$, and let $\l$ be the Lie algebra of $L$.
$L$ admits an IFAS if and only if $[\l, \l] \neq \l$. \\[-3mm]
\end{thm}

Note that Lie algebras of dimension $\leq 5$ satisfying $[\l, \l] = \l$ are exhausted by the following; $\mathfrak{sl}(2,\R)$, $\mathfrak{o}(3)$ and $\mathfrak{sl}(2,\R) \ltimes \R^2$.

\setlength{\tabcolsep}{4pt}

In dimension 6 there is a Lie group $SO(4)$ which does not admit any IFPS. Therefore the minimum dimension of a Lie group which does not admit  any IFPS is 6.


Concerning the existence and non-existence of IFPSs and IFASs on Lie groups, there are several previous works. For example simple Lie groups admitting an IFPS are classified by Agaoka, Urakawa, Elduque (\cite{agaoka1}, \cite{urakawa}, \cite{elduque}). Concerning IFASs, it is known that any 3-step nilpotent Lie group admits an IFAS (Scheuneman \cite{scheuneman}), and any nilpotent Lie group of dimension $\leq 6$ admits an IFAS (Fujiwara \cite{fujiwara}). On the other hand, there are nilpotent Lie groups of dimension $10 \leq n \leq 12$, which do not admit any IFAS (Benoist \cite{benoist}, Burde \cite{burde1}, \cite{burde2}).

To prove the above two Theorems,  we recall a correspondence between the set of IFPSs and the set of certain Lie algebra homomorphisms called (P)-homomorphisms (\cite{agaoka1}) in $\S$2, and recall a classification of real Lie algebras of dimension $\leq 5$ (\cite{mubarakzyanov1}, \cite{mubarakzyanov2}) in $\S$4.  In $\S$3 we give some sufficient conditions for a semidirect sum of Lie algebras to admit an IFAS or IFPS. In $\S$5 and $\S$6 we prove Theorems \ref{IFPS} and \ref{IFAC}
by applying results in $\S$3 to the Lie algebras of dimension $\leq 5$.  
\section{IFPS and (P)-homomorphism}
Let us review some known facts about IFAS and IFPS. First we define a projective structure. Let $\nabla$ and $\nabla'$ be torsion-free affine connections on an $n$-dimensional manifold $M$. Affine connections $\nabla$ and $\nabla'$ are said to be \textit{projectively equivalent} if there exists a 1-form $\phi$ on $M$ such that $\nabla_XY -\nabla'_XY$ $=$ $\phi(X)Y +\phi(Y)X$ for any $X,Y \in \X$, and we express it as $\nabla \sim \nabla'$. A projective equivalence  is an equivalence relation, and the equivalence class $[\nabla]$ containing $\nabla$ is called a \textit{projective structure} on $M$. Now we define IFPS and IFAS on an $n$-dimensional Lie group $L$. We denote by $\nabla_0$ the standard affine connection on $\R^n$ defined by ${\nabla_0}_{\pari}\parj=0$ for $i,j=1,\ldots n$.

\renewcommand{\theenumi}{\roman{enumi}}
\renewcommand{\labelenumi}{\rm (\theenumi)}
\begin{df}\label{def of IFPS}
{\rm $[\nabla]$ {\rm (resp. $\nabla$)} is called an IFPS {\rm (resp. IFAS)} on $L$ if the following two conditions are satisfied. }
\begin{enumerate}
\item {\rm (flatness)} {\rm For each point $p$ of $M$, there exists a neighborhood $U$ of $p$, 
and a diffeomorphism $f$ from $U$ into $\R^n$ such that $f^* \nabla_0 \sim \nabla$ {\rm (resp. 
$f^* \nabla_0 = \nabla$)} on $U$.}
\item {\rm (left invariance)}  {\rm $L_g^*\nabla \sim \nabla$ {\rm (resp. $L_g^*\nabla = \nabla$)} for any $g$ $\in$ $L$. }
\end{enumerate}
\end{df}
Note that an IFAS $\nabla$ naturally induces an IFPS $[\nabla]$.

In the following we define a (P)-homomorphism of the Lie algebra $\l$ of $L$. Then we recall an important theorem concerning the relation between (P)-homomorphisms of $\l$ and IFPSs on $L$, following \cite{agaoka1}. 

Let $\mathfrak{a}_1$ be the one dimensional abelian Lie algebra, and let $\g$ be the special linear algebra $\mathfrak{sl}(\l \oplus \mathfrak{a}_1)$. We fix a basis $e$ of $\mathfrak{a}_1$. 
For any $A \in \mathfrak{sl}(\l \oplus \mathfrak{a}_1)$, there exist $B \in \mathfrak{gl}(\l)$, $u \in \l$ and $\xi \in \l^*$ such that 
\[
\left \{
\begin{array}{ll}
A(x, 0) = (B(x), \ \xi (x)\cdot e) & (\ x \in \l \ ) \\
A(0, e) = (u, \ -\mathrm{tr}B \cdot e) .&  
\end{array}
\right.
\]
Hence we can identify $A$ with 
$
\begin{pmatrix}
B & u \\
\xi & -\mathrm{tr}B
\end{pmatrix}
$.
Then we can decompose $\g$ into $\g_{-1}\oplus \g_{0}\oplus \g_{1}$ defined by

\begin{eqnarray*}
&\g_{-1}&=
\left\{
      \left.
           \begin{pmatrix}
                 0 & u \\
                 0 & 0 
           \end{pmatrix}  
      \right|
          u \in \l
\right\},\\
&\g_0&=
\left\{
     \left.
           \begin{pmatrix}
                 B & 0 \\
                 0 & -\mathrm{tr}B 
           \end{pmatrix}  
      \right|
          B \in \mathfrak{gl}(\l)
\right\},\\
&\g_1&=
\left\{
      \left.
           \begin{pmatrix}
                 0   & 0 \\
                 \xi & 0 
           \end{pmatrix}  
      \right|
          \xi \in \l^*
\right\}.
\end{eqnarray*}
We can see that formal brackets of matrices give this decomposition a graded Lie algebra structure.
The subalgebra $\g_0$ can be identified with the Lie algebra $\mathfrak{gl}(\l)$ under the correspondence; 
\[
\g_0 \ni 
\begin{pmatrix}
B & 0 \\
0 & -\mathrm{tr}B
\end{pmatrix}
\longleftrightarrow B +\mathrm{tr}B\cdot I \in \mathfrak{gl}(\l).
\]
If $f$ is a linear map from $\l$ to $\g$,  then we denote by $f_i$ the $\g_i$-component of $f$. 
\begin{df}[\cite{agaoka1}]
{\rm
Let  $f: \l \to \mathfrak{sl}(\l \oplus \mathfrak{a}_1)$ be a Lie algebra homomorphism. Then $f$ is called a \textit{{\rm(P)}-homomorphism} if $f_{-1}(x) = x$ for any $x \in \l$.

Two (P)-homomorphisms $f$ and $f'$ are said to be \textit{projectively equivalent} if there exists $\xi \in \g_1$ such that $f_0' -f_0 =[\xi,f_{-1}]$. We denote this equivalence relation by $f \sim f'$.
}
\end{df}
If we fix a basis of $\l$, we can easily see that the above definition of (P)-homomorphism is identical with the one in \cite{agaoka1}.

An IFPS on $L$ naturally induces a (P)-homomorphism of $\l$ (\cite{agaoka1}). Furthermore we have the following:
\begin{thm}[\cite{agaoka1}]\label{IFPS and P-hom}
There is a one-to-one correspondence between the set of IFPSs on $L$ and the set of projective equivalence classes of 
{\rm(P)}- homomorphisms of $\l$. 
\end{thm}

\begin{thm}[\cite{agaoka1}]\label{IFAC and P-hom}
There is a one-to-one correspondence between the set of IFASs on $L$ and the set of {\rm(P)}-homomorphisms $f$ of $\l$ satisfying 
$f_1=0$.
\end{thm}
\noindent
\remark
As a representative of an IFPS on $L$ we can take a left invariant projectively flat affine connection(cf. \cite{agaoka1}; Proposition 4.5). 
Furthermore such a left invariant affine 
connection $\nabla$ is uniquely determined for a given (P)-homomorphism $f$ of $\l$ (cf. \cite{agaoka1}; Theorem 3.5). 
In fact we can write explicitly $\nabla$ corresponding to $f$ by $\nabla_xy = f_0(x)y$ for $x, y \in \l$.
If $f$ satisfies the additional condition $f_1 =0$, then from Theorem \ref{IFAC and P-hom}  
$\nabla$ gives an IFAS.\\

From Theorems \ref{IFPS and P-hom} and \ref{IFAC and P-hom} we can decide whether $L$ admits an IFPS (resp.\ IFAS) or not by 
determining whether a (P)-homomorphism  (resp.\ (P)-homomorphism $f$ satisfying $f_{1}=0$) of $\l$ exists or not. If a Lie group $L$ 
admits an IFPS (resp.\ IFAS), we say that the Lie algebra $\l$ admits an IFPS (resp.\ IFAS).

Concerning Theorem \ref{IFAC and P-hom}, we reformulate the condition of (P)-homomorphisms corresponding to IFASs as adapted to 
our calculation.
\begin{lem}\label{lem1}
Let $f: \l \to \mathfrak{sl}(\l \oplus \mathfrak{a}_1)$ be a linear map such that $f_{-1}(x)=x$ and $f_{1}=0$. Then f is a {\rm(P)
}-homomorphism if and only if $f_0 : \l \to \mathfrak{gl}(\l)$ is a Lie algebra homomorphism satisfying the condition
\[\hspace{0cm} f_0(x)y-f_0(y)x =[x, y] \]
for any $x, y \in \l$.
\end{lem}
\pf 
Let $f$ be a linear map $f : \l \to \mathfrak{sl}(\l \oplus \mathfrak{a}_1)$ such that $f_{-1}(x)=x$ and $f_{1}=0$, i.e. $f$ is of the form 
\[
\hspace{0cm}
 f(x) = 
\begin{pmatrix}
A(x) &  x \\
0 & -\mathrm{tr}A(x)
\end{pmatrix}.
\]
Here $A$ is a linear map from $\l$ into $\mathfrak{gl}(\l)$.
Under this situation, $f$ is a (P)-homomorphism  if and only if $f$ is a Lie algebra homomorphism, i.e. the following holds:
\[
(a) \hspace{3pt}
\left \{
\begin{array}{l}
\mbox{$A$ is a Lie algebra homomorphism}, \\
A(x)y+ \mathrm{tr}A(x)\cdot y - (A(y)x + \mathrm{tr}A(y)\cdot x) = [x, y].
\end{array}
\right.
\]
Since $f_0(x)$ equals $A(x) + \mathrm{tr}A(x)\cdot I$, we obtain 
\begin{eqnarray*}
[f_0(x), f_0(y)] &=& [A(x),A(y)],\\
f_0([x, y]) &=& A([x, y]) +\mathrm{tr}A([x, y])\cdot I.
\end{eqnarray*}
Then it is easy to see that the above condition $(a)$ is equivalent to the next condition
\[
\left\{
\begin{array}{l}
\mbox{$f_0$ is a Lie algebra homomorphism}, \\
f_0(x)y - f_0(y)x  = [x, y].
\end{array}
\right.
\]
\hfill $\Box$

Combining Theorem \ref{IFAC and P-hom} with Lemma \ref{lem1}, we obtain the following. This result also can be proved by using the formulas in Chapter 10 of \cite{kobayashi-nomizu}. 
\begin{cor}\label{IFAC and hom}
 Let $L$ be a Lie group. There is a one-to-one correspondence  between the set of IFASs on $L$ and the set of Lie algebra homomorphisms $g$: $\l \to \mathfrak{gl}(\l)$ satisfying
\[ g(x)y-g(y)x=[x, y] \ \ \ \ \  (\ast) \]
for any $x, y \in \l$.
\end{cor}

Concerning Theorem \ref{IFPS and P-hom}, Agaoka \cite{agaoka1} has shown a refined result by using the notion of (N)-homomorphism. 
\begin{df}\rm
Let $f$ be a Lie algebra homomorphism $f : \l \to \mathfrak{sl}(\l \oplus \mathfrak{a}_1)$. Then $f$ is called an \textit{{\rm(N)}-homomorphism} if $f(x)(0, e) = (x, 0)$ for any $x \in \l$. 
\end{df}
In the following we denote simply by $f(x)e = x$ the equality in the above definition. 
By the definition any {\rm(N)}-homomorphism gives a {\rm(P)}-homomorphism, but the converse does not hold.
When we represent $f$ by matrices through an identification of a basis $\{X_1, \cdots, X_n, e \}$ of $\l \oplus \mathfrak{a}_1$ and 
the standard basis $\{e_1, \cdots, e_n, e_{n+1} \}$ of $\R^{n+1}$,
$f : \l \to \mathfrak{sl}(n+1 ,\R)$ is an (N)-homomorphism if and only if $f$ is of the form
\[f(X_i) =
\left(
\begin{array}{c|c}
* & e_i 
\end{array}
\right) \ \ \ (1 \leq i \leq n), \]
where $*$ is an $(n+1) \times n$ matrix.  

It is easy to see that the condition $(\sharp)$ in [17; p. 348] is equivalent to the above definition of (N)-homomorphism. 
\begin{prop}[\cite{agaoka1}]\label{N-hom}
Let $f$ be a {\rm (P)}-homomorphism. Then there exists a unique {\rm (N)}-homomorphism $f'$ such that $f \sim f'$.
\end{prop} 

The next lemma is easy but plays an important role in this paper.
\begin{lem}\label{lem3}
There exists an IFPS on a Lie algebra $\l$ if and only if there exists a Lie algebra homomorphism $f : \l \oplus \mathfrak{a}_1 \to \mathfrak{gl}(\l \oplus \mathfrak{a}_1)$ such that $f(x)e$ $=$ $x$ for any $x \in \l \oplus \mathfrak{a}_1$.
\end{lem}
\pf Suppose that $\l$ admits an IFPS. Then we have the corresponding (N)-homomorphism $g$ of $\l$. We can extend $g$ to $\tilde{g}$ : $\l \oplus \mathfrak{a}_1 \to {\mathfrak{gl}(\l \oplus \mathfrak{a}_1)}$ by setting $\tilde{g}(e) = id$, the identity transformation of $\l \oplus \mathfrak{a}_1$. Then $\tilde{g}$ is a Lie algebra homomorphism and satisfies $\tilde{g}(x)e = x$ for any $x \in \l \oplus \mathfrak{a}_1$. Conversely if $f$ : $\l \oplus \mathfrak{a}_1$ $\to$ $\mathfrak{gl}(\l \oplus \mathfrak{a}_1)$ satisfies the above condition, then $f|_{\l} - \mathrm{tr}(f|_\l)/\mathrm{dim} (\l \oplus \mathfrak{a}_1) \ id$ gives a (P)-homomorphism of $\l$, and hence $\l$ admits an IFPS.  \hfill $\Box$ \\[3mm]
\noindent
\remark
Let $f$ be a homomorphism in Lemma \ref{lem3}. Then from the equality $f([x, y])e=([f(x), f(y)])e$ for $x$, $y$ $\in \l \oplus \mathfrak{a}_1$, we have $f(x)y-f(y)x = [x, y]$.  Thus if $\l$ admits an IFPS, then $\l \oplus \mathfrak{a}_1$ admits an IFAS. 
This gives another proof of Theorem 4.7 in \cite{agaoka1} for the case of Lie groups. Furthermore the equality $f(x)y-f(y)x = [x, y]$ 
induces $f(e)x -f(x)e =0$ for $x \in \l \oplus \mathfrak{a}_1$. From the assumption of $f(x)e=x$, we have $f(e)= id$.

\section{Semidirect sum of Lie algebras admitting an IFAS}
In this section we show that if a Lie algebra $\h$ admits an IFAS and if a nilpotent Lie algebra $\kk$ satisfies some graded condition, then the semidirect sum $\h \ltimes \kk$ also admits an IFAS (Proposition \ref{solv and IFAC}).  By using this proposition,  we can show that most Lie algebras of dimension $\leq 5$  admit an IFAS. 

In the following we denote by $Z(\kk)$ the center of $\kk$, and by $\h \ltimes \kk$ a semidirect sum of  $\kk$ by $\h$. Let $g : \l \to \mathfrak{gl}(\l)$ be a Lie algebra homomorphism. Then we express $g(x)y$ simply as $x\cd y$ for $x$, $y$ $\in \l$.
\begin{prop}\label{nil and IFAC}
Let $\kk$ be a Lie algebra which has a direct sum decomposition $ \displaystyle \kk = \oplus_{i \geq 1} \kk_i \oplus Z'(\kk)$ as vector spaces such that $[\kk_i, \kk_j]$ $\subset \kk_{i+j} \oplus Z'(\kk)$, and $Z'(\kk)$ is a subspace of $Z(\kk)$. Then $\kk$ admits an IFAS.
\end{prop}
\pf
Let us define a linear map $g : \kk \to \mathfrak{gl}(\kk)$ by
\[
x \cd y = 
\left \{
\begin{array}{cl}
\frac{j}{i+j}[x,y] & x \in \kk_i, \ y \in \kk_j \\
0 & x \in \kk, \ y \in Z'(\kk) \\
0 & x \in Z'(\kk), \ y \in \kk.
\end{array}
\right.
\]
Then g satisfies the condition $(\ast)$ in Corollary \ref{IFAC and hom} : $x \cd y$ $-$ $y \cd x$ $= [x, y]$  for $x$, $y$ $\in \kk$.

We show that $g$ is a Lie algebra homomorphism, i.e. $x \cd (y \cd z) - y \cd (x \cd z) - [x,y] \cd z = 0$. We have to check three cases, (1) $x \in \kk_i$, $y \in \kk_j$, $z \in \kk_k$, \ (2) $x \in \kk$, $y \in \kk$, $z \in Z'(\kk)$, \  (3) $x \in \kk$, $y \in Z'(\kk)$, $z \in \kk$. 
We can verify (2) and (3) easily, thus we check only the case (1).  

(1) \  Because of the condition $[\kk_i, \kk_j]$ $\subset \kk_{i+j} \oplus Z'(\kk)$, we can decompose $[y, z]$ into $[y, z]_a + [y, z]_b$ $\in$ $\kk_{j+k} \oplus Z'(\kk)$.  Note that $[x, [y, z]] = [x, [y, z]_a]$. Then
\begin{eqnarray*}
\hspace{5mm} && x \cd (y \cd z) - y \cd (x \cd z) - [x,y] \cd z \\
&=& x \cd \frac{k}{j+k} ([y,z]_a + [y,z]_b)- y \cd \frac{k}{i+k}( [x,z]_a + [x,z]_b) \\
&& {} - ([x,y]_a + [x,y]_b) \cd z \\
&=& \frac{j+k}{i+j+k}\frac{k}{j+k}[x,[y,z]_a] - \frac{i+k}{j+i+k} \frac{k}{i+k}[y,[x,z]_a] -  \frac{k}{i+j+k}[[x,y]_a,z] \\
&=& \frac{k}{i+j+k} [x,[y,z]] - \frac{k}{i+j+k} [y,[x,z]] - \frac{k}{i+j+k} [[x,y],z] \\
&=& 0.
\end{eqnarray*}

\noindent
Hence $\kk$ admits an IFAS. \hfill $\Box$\\[3mm]
\remark
A Lie algebra $\kk$ which satisfies the condition in Proposition \ref{nil and IFAC} is necessarily  nilpotent.
By Proposition \ref{nil and IFAC} we can verify that any three-step nilpotent Lie algebra $\kk$ admits an IFAS (Scheuneman \cite{scheuneman}) as follows. We define an ideal $C^i\kk$ of $\kk$ by $C^i\kk$ $:=$ $[\kk, C^{i-1}\kk]$ ($i=2, 3, \cdots$), $C^1\kk$ $=$ $\kk$ ($i=1$). Let us denote complementary subspaces $C^{i+1}\kk$ in $C^{i}\kk$ by $\kk_i$, i.e., $C^{i}\kk$ $=$ $\kk_i \oplus C^{i+1}\kk$. Now we take $Z'(\kk) = \kk_3$. Then the decomposition $\kk$ $=$ $\kk_1 \oplus \kk_2 \oplus Z'(\kk)$ satisfies the condition in Proposition \ref{nil and IFAC}. 

Moreover this proposition  improves the following result of Burde \cite{burde3}: If a Lie algebra $\kk$ is graded by positive integers, then $\kk$ admits an IFAS. We can verify that a Lie algebra $ \displaystyle \kk = \oplus_{i \geq 1} \kk_i \oplus Z'(\kk)$ satisfies the condition $[\kk_i, \kk_j] \subset \kk_{i+j} \oplus Z'(\kk)$ if and only if $\kk \big/ Z'(\kk)$ is graded by positive integers.
\begin{prop}\label{solv and IFAC}
Let $\h$ be a Lie algebra which admits an IFAS. Let $\kk$ be a Lie algebra which has a decomposition $\kk = \bigoplus_{i \geq 1} \kk_i \oplus Z'(\kk)$ such that  $[\kk_i, \kk _j]$ $\subset$ $\kk_{i+j} \oplus Z'(\kk)$, and $Z'(\kk) \subset Z(\kk)$.  Then any semidirect sum $\l = \h \ltimes \kk$ satisfying $[\h, \kk_i] \subset \kk_i \oplus Z'(\kk)$ and $[\h, Z'(\kk)] \subset Z'(\kk)$ admits an IFAS.
\end{prop}
\pf 
A given IFAS on $\h$ corresponds to a Lie algebra homomorphism $g : \h \to \mathfrak{gl}(\h)$ satisfying the condition $(\ast)$ in Corollary \ref{IFAC and hom} : $x \cd y$ $-$ $y \cd x$ $=$ $[x ,y]$ for $x$, $y$ $\in \h$.  Let $h :\kk \to \mathfrak{gl}(\kk)$ be the homomorphism constructed in the proof of Proposition \ref{nil and IFAC}. By using g and h we define a linear map $f : \l \to \mathfrak{gl}(\l)$ by the following 
\[
x \cd y = 
\left \{
\begin{array}{cl}
g(x)y & x \in \h, \ y \in \h \\
h(x)y & x \in \kk, \ y \in \kk \\
\left[x, y \right] & x \in \h, \ y \in \kk \\
0 & x \in \kk, \ y \in \h. \\
\end{array}
\right.
\]
\noindent
Then $f$ satisfies the condition $(\ast)$: $x\cd y$ $-$ $y \cd x$ $=$ $[x, y]$ for $x$, $y$ $\in \l$.

We show that $f$ is a Lie algebra homomorphism, i.e. $x \cd (y \cd z)$ $-$ $y \cd (x \cd z)$ $-$ $[x,y] \cd z$ $=$ $0$. Then it is enough to check six cases;  
(1) $x \in \h$, $y \in \h$, $z \in \kk$, 
(2) $x \in \h$, $y \in \kk$, $z \in \h$,
(3) $x \in \h$, $y \in \kk_i$, $z \in \kk_j$,
(4) $x \in \h$, $y \in \kk$, $z \in Z'(\kk)$,
(5) $x \in \h$, $y \in Z'(\kk)$, $z \in \kk$,
(6) $x \in \kk$, $y \in \kk$, $z \in \h$.

Here we verify only (3). 
\begin{eqnarray*}
(3)  && x \cd (y \cd z) - y \cd (x \cd z) - [x,y] \cd z \\
&=& x \cd (\frac{j}{i+j} [y, z]) - y \cd [x, z] - [x,y] \cd z \\
&=& \frac{j}{i+j}[x,[y,z]] - \frac{j}{i+j}[y,[x,z]] - \frac{j}{i+j}[[x,y],z]\\
&=& 0.\\
\end{eqnarray*}
The other cases also can be proved by the straightforward computation. 
Thus $\l$ admits an IFAS. \hfill $\Box$\\[2mm]

\remark
There are several  preceding results related to Proposition \ref{solv and IFAC} as follows: If $H$ is a solvable affinely flat Lie group and $K$ = $\R^n$, then the semidirect product $L = H \ltimes K$ admits an IFAS (Mizuhara \cite{mizuhara}). If $\h$ is affinely flat and $\kk$ is an abelian Lie algebra, then $\h \ltimes \kk$ admits an IFAS (Burde \cite{burde2}). Moreover if $\kk$ is a 2-step nilpotent Lie algebra, then $\h \ltimes \kk$  admits an IFAS (Fujiwara \cite{fujiwara}).  Suppose that $\h$ is abelian, and $\kk$ is graded by positive integers. If $\l = \h \ltimes \kk$ satisfies $[\h, \kk_i] \subset \kk_i$, then $\l$ admits an IFAS (Yamaguchi \cite{yamaguchi}). As can be easily seen, Proposition \ref{solv and IFAC} is an extension of these results.\\

Next we consider a direct sum of Lie algebras admitting an IFAS or IFPS. It is well known that from two affine manifolds $(M_i, \nabla_i)$ ($i=1, 2$) we can construct the product affine manifold ($M_1 \times M_2, \nabla_1 \times \nabla_2$). Moreover if ($M_i, \nabla_i$) are torsion-free flat affine manifolds, the product is also torsion-free flat affine manifold. Thus the next proposition is immediately obtained from this fact.  However we prove this from the viewpoint of representation.  
\begin{prop}\label{sum of IFAC1}
If each Lie algebra $\mathfrak{l}_1$ and $\mathfrak{l}_2$ admits an IFAS, then the direct sum $\mathfrak{l}_1 \oplus \mathfrak{l}_2$ 
also admits an IFAS.
\end{prop}
\pf 
From Corollary \ref{IFAC and hom} there exist Lie algebra homomorphisms $g_i : \mathfrak{l}_i\rightarrow \mathfrak{gl}(\mathfrak{l}_i)$ $(i=1,2)$ such that $g_i(x)y-g_i(y)x=[x,y]$ ($x$, $y$ $\in \mathfrak{l}_i$). Let us define a map $g:\mathfrak{l}_1\oplus \mathfrak{l}_2 \rightarrow \mathfrak{gl}(\mathfrak{l}_1)\oplus \mathfrak{gl}(\mathfrak{l}_2)$ $\subset$ $\mathfrak{gl}(\l_1 \oplus \l_2)$ by $g(x_1,x_2)$ $=$ $(g_1(x_1),g_2(x_2))$. Then it follows easily that $g$ is a Lie algebra homomorphism and satisfies $g((x_1,x_2))(y_1,y_2) -g((y_1,y_2))(x_1,x_2)$ $=$ $[(x_1,x_2),(y_1,y_2)]$. Hence by Corollary \ref{IFAC and hom} $\mathfrak{l}_1\oplus \mathfrak{l}_2$ admits an IFAS. \\ 
\ \ \hfill $\Box$\\
  
Even if each $\mathfrak{l}_1$ and $\mathfrak{l}_2$ admits an IFPS, the direct sum  $\mathfrak{l}_1\oplus \mathfrak{l}_2$ does not necessarily admit an IFPS. Indeed there exists a counterexample. Although $\mathfrak{o}(3)$ admits an IFPS (cf. $\S$5), the direct sum $\mathfrak{o}(3) \oplus \mathfrak{o}(3)$ does not admit any IFPS. Because if $\mathfrak{o}(3) \oplus \mathfrak{o}(3)$ admits an IFPS, then the Lie group $S^3 \times S^3$ admits an IFPS. But this contradicts the following result by Kobayashi-Nagano \cite{kobayashi-nagano}: A compact simply connected manifold with flat projective structure is diffeomorphic to $S^n$.  Thus $\mathfrak{o}(3) \oplus \mathfrak{o}(3)$ does not admit any IFPS. However we have the following.

\begin{prop}\label{direct sum of IFPS}
Suppose that each $\l_1$ and $\l_2$ admits an IFPS. Then $\l_1 \oplus \l_2 \oplus \mathfrak{a}_1$ admits an IFPS.
\end{prop}
\pf 
Let $\mathfrak{a}_1$, $\mathfrak{a}_2$ be the one dimensional abelian Lie algebras, and let $e_1$ and $e_2$ be bases of $\mathfrak{a}_1$ and $\mathfrak{a}_2$. We define two subalgebras $\l'$ and $\mathfrak{a}$ of $\l_1 \oplus 
\mathfrak{a}_1 \oplus \l_2 \oplus \mathfrak{a}_2$ by 
\[
\left \{
\begin{array}{l}
\l' = \{ (x, k e_1, y, -k e_2) \ |\  x \in \l_1, y \in \l_2, k \in \R \} \\
\mathfrak{a} = \{ (0, k e_1, 0, ke_2) \ |\  k \in \R \}.
\end{array}
\right.
\]
Here we have obviously $\l' \oplus \mathfrak{a}$ = $\l_1 \oplus \mathfrak{a}_1 \oplus \l_2 \oplus \mathfrak{a}_2$.  
From Lemma \ref{lem3} we have homomorphisms $f_i$ $:$ $\l_i \oplus \mathfrak{a}_i  \to \mathfrak{gl}(\l_i \oplus \mathfrak{a}_i)$  such that $f_i(x_i) e_i = x_i$ for any $x_i \in \l_i \oplus \mathfrak{a}_i$ ($i =1, 2$).
We show that $f = f_1 \oplus f_2$ : $(\l_1 \oplus \mathfrak{a}_1) \oplus (\l_2 \oplus \mathfrak{a}_2)$ $\to$ $\mathfrak{gl}(\l_1 \oplus \mathfrak{a}_1) \oplus \mathfrak{gl}(\l_2 \oplus \mathfrak{a}_2)$ $\subset$ $\mathfrak{gl}((\l_1 \oplus \mathfrak{a}_1) \oplus (\l_2 \oplus \mathfrak{a}_2))$ gives an IFPS on $\l'$. 
If we choose $(0,  e_1, 0, e_2)$ as a basis of $\mathfrak{a}$, then $f : \l' \oplus \mathfrak{a} \to \mathfrak{gl}(\l' \oplus \mathfrak{a})$ satisfies the condition in Lemma \ref{lem3}. Hence $\l'$ admits an IFPS. Since $\l'$ is isomorphic to $\l_1 \oplus \l_2 \oplus \mathfrak{a}_1$, it follows that $\l_1 \oplus \l_2 \oplus \mathfrak{a}_1$ admits an IFPS.
\hfill $\Box$ \\

\section{Classification of Lie algebras of dimension $\leq$ 5}
Mubarakzyanov has classified real Lie algebras of dimension $\leq$ 5 in \cite{mubarakzyanov1}, \cite{mubarakzyanov2}. This result is restated in Patera-Sharp-Winternitz-Zassenhaus \cite{patera}, which we recall below, and we use this classification in the proof of Theorems \ref{IFPS} and \ref{IFAC}. 

Note that in the following list Lie algebras are denoted by the same symbol $A_{n,i}$ as \cite{patera}, and    
the direct sums of Lie algebras are also omitted. Hence 
the Lie algebras which are not direct sums of  lower dimensional Lie algebras (or indecomposable Lie algebras, for short) are classified.  We denote by $\{X_1,\ldots, X_n\}$  a basis of an $n$-dimensional Lie algebra, and $[X_i,X_j]$ ($i < j$) its bracket.  We omit the bracket $[X_i,X_j]$ when $[X_i,X_j]=0$. The basis $\{X_1,\ldots, X_n\}$ is different from the one $\{e_1, \ldots, e_n\}$ in \cite{patera}, and  
the relation is given by  $\{X_1,X_2, \ldots, X_n\} = \{-e_n, e_{n-1}, \ldots, e_1\}$ except for $A_{1,1}$, $A_{2,1}$ and  $A_{3,8}$ ($\sll(2, R)$). For example $\{X_1,X_2, X_3, X_4\} = \{-e_4, e_3, e_2, e_1\}$. 
\\[5mm]
\newpage
\noindent
dimension 1 \quad $A_{1,1}$ \ abelian Lie algebra\\
dimension 2 \quad $A_{2,1}$ \ $[X_1,X_2]=X_2$ \ \ (solvable)\\
\noindent
dimension 3
\vspace{-0.1cm}
\begin{table}[htbp]
\begin{tabular}{llll}
$A_{3,1}$ & $[X_1, X_2]=X_3,$ \\[1mm]
$A_{3,2}$ & $[X_1, X_2]=X_2 + X_3$, & $[X_1, X_3]=X_3$, \\[1mm]
$A_{3,3}, A_{3,4}, A_{3,5}^a$ & $[X_1, X_2]=a X_2,$ & $[X_1, X_3]=X_3$ & $(0<|a|\leq  1),$ \\[1mm]
$A_{3,6}, A_{3,7}^a$ & $[X_1, X_2]= a X_2+X_3,$ & $[X_1, X_3]=-X_2 + a X_3$ & $(a \geq 0),$ \\[1mm]
$A_{3,8} \ (\mathfrak{sl}(2,\R))$ & $[X_1, X_2]= X_2,$ & $[X_1, X_3]= -X_3,$ & $[X_2, X_3]=X_1,$ \\[1mm]
$A_{3,9} \ (\mathfrak{0}(3))$ & $[X_1, X_2]=X_3,$ & $[X_2, X_3]=X_1,$ & $[X_3, X_1]=X_2.$
\end{tabular}
\end{table}

\noindent
$A_{3,1}$ is  nilpotent, and $A_{3,2} \sim A_{3,7}^a$ are solvable.\\[2mm]
dimension 4
\vspace{-0.1cm}
\setlength{\tabcolsep}{4pt}
\begin{table}[htbp]
\begin{tabularx}{\linewidth}{lllll}
$A_{4,1}$ & $[X_1,X_2]\!=\!X_3,$ & $[X_1,X_3]\!=\!X_4,$ \\[1mm]
$A_{4,2}^a$ & $[X_1,X_2]\!=\!X_2 + X_3$, & $[X_1,X_3]\!=\!X_3$, & $[X_1,X_4]\!=\! aX_4$,\\
&&& $(a \neq 0)$, \\[1mm] 
$A_{4,3}$ & $[X_1,X_2]\!=\!X_3,$  & $[X_1,X_4]\!=\!X_4$, \\[1mm]
$A_{4,4}$ & $[X_1,X_2]\!=\!X_2+X_3,$ & $[X_1,X_3]\!=\!X_3 +X_4,$ & $[X_1,X_4]\!=\!X_4,$ \\[1mm]
$A_{4,5}^{ab}$ & $[X_1,X_2]\!=\!b X_2,$ & $[X_1,X_3]\!=\!aX_3,$ & $[X_1,X_4]\!=\!X_4$ &\hspace{6cm} \\
&&&$(ab \neq 0, -1\!\leq \! a \! \leq \! b \! \leq \!1),$ \\[1mm]
$A_{4,6}^{ab}$ & $[X_1,X_2]\!=\!2 X_2 + X_3,$ & $[X_1,X_3]\!=\!-X_2 +b X_3$, & $[X_1,X_4]\!=\!a X_4$ \\
&&& $(a\neq 0, b\geq 0),$ \\[1mm]
$A_{4,7}$ & $[X_1,X_2]\!=\!X_2+X_3,$ & $[X_1,X_3]\!=\!X_3,$ & $[X_1,X_4]\!=\!2X_4,$ \\ & 
$[X_2,X_3]=- X_4,$ \\[1mm]
$A_{4,8}, A_{4,9}^b$ & $[X_1,X_2]\!=\!b X_2,$ & $[X_1,X_3]\!=\!X_3,$ & $[X_1,X_4]\!=\!(1\!+\!b)X_4$, \\& $[X_2,X_3]=- X_4$ && $(-1\leq b \leq 
1),$\\[1mm]
$A_{4,10}, A_{4,11}^a$ & $[X_1,X_2]\!=\!aX_2+X_3,$ & $[X_1,X_3]\!=\!-X_2+aX_3$, & $[X_1,X_4]\!=\!2aX_4,$\\ & $[X_2,X_3]=- X_4$ &&$(a\geq 0),$\\
[1mm]
$A_{4,12}$ & $[X_1,X_3]\!=\!X_4,$ & $[X_1,X_4]\!=\!-X_3,$ & $[X_2,X_3]\!=\!-X_3,$ \\& $[X_2,X_4]=-X_4.$ 
\end{tabularx}
\end{table}

\vspace{-0.2cm}
\noindent
$A_{4,1}$ is nilpotent, and $A_{4,2}^a \sim A_{4,12}$ are solvable.\\

In the above list some different Lie algebras are defined together by using parameters $a$, $b$.   
For example $A_{3,3}$, $A_{3,4}$ and $A_{3,5}^{a}$ is distinguished 
by  $a = 1$, $a =-1$ and $0<|a|\leq  1$ in the above definition. 
In this case $A_{3,5}^{a}$ denotes a family of Lie algebras. 
In dimension 5, indecomposable Lie algebras consist of 40 families $A_{5,1} \sim A_{5,40}$ (see \cite{patera} for their definitions). Among 
them six Lie algebras $A_{5,1} \sim A_{5,6}$ are nilpotent, and $A_{5,40}$ is isomorphic to the semidirect sum $\mathfrak{sl}(2, \R) 
\ltimes \R^2$ which is perfect. The remaining Lie algebras $A_{5,7} \sim A_{5,39}$ are all solvable.
\\

\section{Proof of Theorems \ref{IFPS} and \ref{IFAC} (indecomposable case)}
In this and next sections, we prove Theorems \ref{IFPS} and \ref{IFAC}. If a Lie algebra $\l$ admits an IFAS, then $\l$ admits an IFPS. So we study first the existence of IFAS and next the existence of IFPS.
In this section we classify indecomposable Lie algebras of dimension $\leq 5$ admitting IFASs (IFPSs). For the remaining case, i.e. direct sums of Lie algebras, we study them in the next section. \\

Let $\l$ be an indecomposable Lie algebra of dimension $\leq 5$. Then we may suppose that $\l$ is isomorphic to one of the Lie algebras in the list of $\S 4$ if $\dim \l \leq 4$, or the list in 
\cite{patera} if $\dim \l =  5$. First let us examine IFASs.
If $\l$ is perfect, i.e. $[\l, \l] = \l$, then $\l$ admits no IFAS (Helmstetter \cite{helmstetter}).
In the following we assume $\l$ is not perfect. Then we can show that $\l$ has a semidirect sum structure $\l = \h \ltimes \kk$ 
satisfying the condition in Proposition \ref{solv and IFAC}. 
Each semidirect sum structure is given in the following table.\\

\setcounter{table}{-1}
\hspace{3cm} {\sc Table}: \ Non-perfect Lie algebras of dimension $\leq 5$ and \quad \quad\\
\hspace{3.5cm} their semidirect sum structures
\vspace{2mm}
\begin{center}
\begin{tabular}{| c | l | l | l | l | l | l |}
\hline
dim & \multicolumn{1}{c|}{Lie algebra} & \multicolumn{1}{c|}{$\h$} & \multicolumn{1}{c|}{$\kk_1$} & \multicolumn{1}{c|}{$\kk_2$}  & \multicolumn{1}{c|}{$\kk_3$}\\
\hline
1 & $A_{1,1}$ & & $X_1$ \hfill  & & \\
\hline
2 & $A_{2,1}$ & $X_1$ \quad & $X_2$ \quad & &\\
\hline
3 & $A_{3,1} \sim A_{3,7}^a$ & $X_1$ & $X_2, X_3$ & &\\
\hline
4 & $A_{4,1} \sim A_{4,6}^{ab}$ & $X_1$ & $X_2, X_3, X_4$ & &\\
\cline{2-6}
& $A_{4,7} \sim A_{4,11}^a$ & $X_1$ & $X_2, X_3$ & $X_4$ & \\
\cline{2-6}
 & $A_{4,12}$ & $X_1, X_2$ & $X_3, X_4$ & & \\
 \hline
\cline{2-6} 
5 & $A_{5,1}$, $A_{5,2}$, $A_{5,7} \sim A_{5,18}$ & $X_1$ & $X_2, X_3, X_4, X_5$ & & \\
\cline{2-6}
 & $A_{5,3}, A_{5,4}$, $A_{5,30}, A_{5,31}$ & $X_1$ & $X_2, X_3$ & $X_4$ & $X_5$ \\
\cline{2-6}
& $A_{5,5}$, $A_{5,6}$, $A_{5,19} \sim A_{5,29}$ & $X_1$ & $X_2, X_3, X_4$ & $X_5$ & \\
\cline{2-6}
 & $A_{5,32} \sim A_{5,37}$ & $X_1, X_2$ & $X_3, X_4$ & $X_5$ & \\
\cline{2-6}
 & $A_{5,38}$, $A_{5,39}$ & $X_1, X_2, X_3$ & $X_4, X_5$ & & \\
 \hline
\end{tabular}
\end{center}
\vspace{4mm}

We have two comments about the table. First we note that the above basis $\{X_1, \ldots, X_n\}$ of $A_{n,i}$  is different from the one 
$\{e_1, \ldots, e_n\}$ in \cite{patera}, and there is the common relation for $n \geq 3$: $\{X_1, X_2, \ldots  X_n\} =  \{-e_n, e_{n-1}, \ldots, e_1\}$. 
Thus for $n = 5$ we have $\{X_1,X_2, X_3, X_4, X_5\} = \{-e_5, e_4, e_3, e_2, e_1\}$. 
Secondly 
the last non-zero $\kk_i$ of each row stands for a subspace $Z'(\kk)$ of the center $Z(\kk)$.
For instance $A_{5,3} = \h \oplus \kk_1 \oplus \kk_2 \oplus Z'(\kk)$, where $Z'(\kk) =\langle X_5 \rangle$.

Here let us examine the decomposition in the table by taking up the Heisenberg Lie algebra $A_{3,1}$. In this case $\h = \langle X_1 
\rangle$ is abelian, so $\h$ admits an IFAS. The ideal $\kk$ $=$ $Z'(\kk)$ $=$ $\langle X_2, X_3 \rangle$ is also abelian, and 
obviously we have $[\h, Z'(\kk)] \subset Z'(\kk)$. Hence $A_{3,1}$ admits an IFAS. Next we take up $A_{5,39}$ defined below.\\[1mm]
\setlength{\tabcolsep}{4pt}
\begin{tabular}{lllll}
$A_{5,39}$ & $[X_1,X_2]\!=\! -X_3$, &  & $[X_1,X_4]\!=\!-X_5$, & $[X_1, X_5] \! = \! X_4$, \\[1mm]
& &  & $[X_2,X_4]\!=\! -X_4$, & $[X_2,X_5]\!=\! -X_5$.
\end{tabular}
\\[3mm]
Then $\h = \langle X_1, X_2, X_3 \rangle$ is isomorphic to $A_{3,1}$, so $\h$ admits an IFAS. The ideal $\kk$ $=$ $Z'(\kk)$ $=$ 
$\langle X_4, X_5 \rangle$ is abelian, and satisfies $[\h, Z'(\kk)] \subset Z'(\kk)$. Hence $A_{5, 39}$ is an extension of $A_{3,1}$ by a 
semidirect sum satisfying the condition in Proposition \ref{solv and IFAC}. In this way for all Lie algebras in the list of $\S$4 and  
\cite{patera}, we can check that they have decompositions in the table and satisfy the condition in Proposition \ref{solv and IFAC}. Hence 
any non-perfect indecomposable Lie algebra of dimension $\leq 5$ admits an IFAS. \\[2mm]
\remark
The existence of IFAS on the above non-perfect Lie algebras can be showed by using the results of 
\cite{fujiwara} and \cite{yamaguchi} which we 
stated in the remark after Proposition \ref{solv and IFAC}. 
But without Proposition \ref{solv and IFAC} it becomes harder to verify  
the existence of IFAS, since the number of case by case 
examinations increases. 
Moreover consider the following 6-dimensional Lie algebra 
$n6\underline{\ }20\underline{\ }1$ (cf. \cite{komrakov}; page 156): \\[2mm]
\begin{tabular}{lllll}
\ $n6\underline{\ }20\underline{\ }1$ & $[e_1, e_3]\!=\! e_3$, & $[e_1, e_5]\!=\! e_2$, & $[e_1, e_6] \! = \! e_5$, \\[1mm]
& $[e_4, e_6]\!=\! e_2$, & $[e_5,e_6]\!=\! e_4$.
\end{tabular}\\[2mm]
We can show the existence of an IFAS on $n6\underline{\ }20\underline{\ }1$ by using Proposition \ref{solv and IFAC} as follows:
$\h = \langle e_1, e_3 \rangle$, $\kk_1 = \langle e_5, e_6 \rangle$, $\kk_2 = \langle e_4 \rangle$, $Z'(\kk) = \langle e_2 \rangle$.
But it seems impossible to show the existence of an IFAS without using Proposition \ref{solv and IFAC}.\\

Next we shall study the existence and non-existence of IFPSs for the remaining perfect indecomposable Lie algebras of dimension $\leq 5$. 
Such Lie algebras are exhausted by $\mathfrak{sl}(2,\R)$, $\mathfrak{o}(3)$ and $A_{5,40}$. 
Agaoka \cite{agaoka1} has shown that each $\mathfrak{sl}(2,\R)$ and $\mathfrak{o}(3)$ admits an IFPS. The corresponding 
(N)-homomorphisms
are given by the following:
\\[2mm]
\renewcommand{\arraystretch}{1.1}
\noindent
{\arraycolsep = 1.4mm
$\mathfrak{sl}(2,\R)$:\\[3mm]
{\scriptsize
$f(X_1) \!=\!
\left(
\begin{array}{cccc}
0& 0 & 0 & 1 \\
0& \dfrac{1}{2}& 0& 0 \\
0& 0& -\dfrac{1}{2}& 0 \\
\dfrac{1}{4}& 0 & 0 & 0 
\end{array}
\right)$, \ 
$f(X_2) \!=\!
\left(
\begin{array}{cccc}
0& 0 & \dfrac{1}{2} & 0 \\
-\dfrac{1}{2}&0& 0& 1\\
0& 0 & 0 & 0\\
0 & 0 & \dfrac{1}{4} & 0 
\end{array}
\right)$, \ 
$f(X_3) \!=\!
\left(
\begin{array}{cccc}
0&-\dfrac{1}{2} & 0 & 0 \\
0& 0& 0& 0\\
\dfrac{1}{2} & 0 & 0 & 1\\
0 & \dfrac{1}{4} & 0 & 0 
\end{array}
\right)$.
}
\\[0.5cm]
$\mathfrak{o}(3)$:\\
{\scriptsize
$f(X_1) \!=\!
\left(
\begin{array}{cccc}
0& 0 & 0 & 1 \\
0& 0& -\dfrac{1}{2}& 0 \\
0& \dfrac{1}{2}& 0& 0 \\
-\dfrac{1}{4}& 0 & 0 & 0 
\end{array}
\right), \ 
f(X_2) \!=\!
\left(
\begin{array}{cccc}
0& 0 & \dfrac{1}{2} & 0 \\
0&0& 0& 1\\
-\dfrac{1}{2} & 0 & 0 & 0\\
0 &-\dfrac{1}{4} & 0 & 0 
\end{array}
\right), \ 
f(X_3) \!=\!
\left(
\begin{array}{cccc}
0&-\dfrac{1}{2} &0 & 0 \\
\dfrac{1}{2}& 0& 0& 0\\
0&0 & 0 & 1\\
0 & 0 &-\dfrac{1}{4} & 0 
\end{array}
\right)$.
}
\\[4mm]
Here $\{X_1, X_2, X_3\}$ is a basis defined in $\S 4$. 
For each case 
left invariant projectively flat affine connection $\nabla$ corresponding to $f$  
is given by $\nabla_xy =f_0(x)y= 1/2 [x, y]$.\\[-2mm]

Finally we consider the remaining Lie algebra $A_{5,40}$. In \cite{patera} $A_{5,40}$ is defined by the following: \\[3mm]
\begin{tabular}{lllll}
$[e_1,e_2]\!=\!2e_1$, & $[e_1,e_3]\!=\!-e_2$, & $[e_2,e_3]\!=\!2e_3$, \\
$[e_1, e_4] \! = \!  e_5$, & $[e_2, e_4]\!=\!e_4$, & $[e_2,e_5]\!=\!-e_5$, & $[e_3,e_5]\!=\!e_4.$
\end{tabular}
\\[3mm]
\noindent
The Lie algebra $A_{5,40}$ is isomorphic to the matrix Lie algebra defined by 
\[\left\{
\left(
\begin{array}{c|c}
A & v \\
\hline
0 & 0
\end{array}
\right)
| A \in \sll(2, \R), v \in \R^2
\right\}\]
through the correspondence between bases
$\{e_1, e_2, e_3, e_4, e_5\}$ and 
\[\left\{
\footnotesize 
\begin{pmatrix}
0 & 0 & 0 \\
1 & 0 & 0\\
0 & 0 &0
\end{pmatrix}, 
\begin{pmatrix}
1 & 0 & 0 \\
0 & -1 & 0\\
0 & 0 &0
\end{pmatrix}, 
\begin{pmatrix}
0 & 1 & 0 \\
0 & 0 & 0\\
0 & 0 &0
\end{pmatrix}, 
\begin{pmatrix}
0 & 0 & 1 \\
0 & 0 & 0\\
0 & 0 &0
\end{pmatrix},  
\begin{pmatrix}
0 & 0 & 0 \\
0 & 0 & 1\\
0 & 0 &0
\end{pmatrix}
\right\}.
\]
Thus the Lie algebra $A_{5,40}$ has a semidirect sum structure $\sll(2, \R) \ltimes \R^2$. 

From the above definition of $\sll(2, \R) \ltimes \R^2$ by matrices, we obtain the standard representation $\imath: \sll(2, \R) \ltimes 
\R^2 \to \gl(3, \R)$. 
From $\imath$ we define a Lie algebra representation $\imath \oplus \imath: \sll(2, \R) \ltimes \R^2 \to \gl(6, \R)$ by
$\imath \oplus \imath (x) (\xi, \eta) = (\imath(x)\xi, \imath(x)\eta)$ for $\xi, \eta \in \R^3$. 
Here we identify $\R^6$ with the set of $3 \times 2$ real matrices. 

The representation $\imath \oplus \imath$ does not give an
(N)-homomorphism, however in this case its contragredient representation $(\imath \oplus \imath)^*$ does. We denote it by $f$, 
i.e. $f(x)(\xi, \eta) = (-{}^t\imath(x) \xi, -{}^t\imath(x) \eta)$ for $\xi, \eta \in \R^3$.
We denote  by $\mathfrak{a}_1$ a one-dimensional vector space, and denote by $e$ its basis. Then we extend $f$ 
to a representation $\tilde{f}: (\sll(2, \R) \ltimes \R^2) \oplus \mathfrak{a}_1 \to \gl(6, \R)$ by mapping $e$ into  
the identity map of 
$\R^6$.
When we set 
$3 \times2$ real matrix \[v =
\begin{pmatrix}
I_2 \\
\hline
0
\end{pmatrix},\] 
$v$ satisfies $\tilde{f}((\sll(2, \R) \ltimes \R^2) \oplus \mathfrak{a}_1)v = \R^6$. Because of this equality, we can identify 
$(\sll(2, \R) \ltimes \R^2) \oplus \mathfrak{a}_1$ with $\R^6$ by the correspondence $x \leftrightarrow \tilde{f}(x)v$.
Under this identification we can regard $f$ as a representation 
\[f: \sll(2, \R) \ltimes \R^2 \to \gl((\sll(2, \R) \ltimes \R^2) \oplus \mathfrak{a}_1).\] 
Then
$f$ satisfies $f(x)e =x$ and $\mathrm{tr}f(x) =0$ for $x \in \sll(2, \R) \ltimes \R^2$, which implies that $f$ is an (N)-homomorphism.

The corresponding left invariant projectively flat affine connection is given by $\nabla_xy = f_0(x)y$ $=$ $-yx + \tr(yx)/2 \ I_2$ for 
$x, y \in \sll(2, \R)$,  $\nabla_ux = -xu$ for $x \in \sll(2, \R)$, $u \in \R^2$, and $\nabla = 0$ for other cases.  
Thus this connection does not coincide with the connection $\nabla_xy = 1/2 [x, y]$.
Note that the above construction can be naturally extended to $\sll(n, \R) \ltimes \R^n$ $(n \geq 2)$, namely $\sll(n, \R) \ltimes \R^n$
 admits an IFPS.

From the above considerations, we complete the proof of Theorems \ref{IFPS} and \ref{IFAC} for indecomposable Lie algebras.
\section{Proof of Theorems \ref{IFPS} and \ref{IFAC} (decomposable case)}
In this section we shall show that every decomposable Lie algebra of
dimension $\leq 5$ admits an IFAS. 
For this purpose we prove Proposition
\ref{sum of IFAC2}.  First we introduce the notion of reducible IFAS.
\renewcommand{\labelenumi}{\rm(\theenumi)}
\begin{df}
\rm
Suppose that $\nabla$ is an IFAS on a Lie algebra $\h$, and $g: \h \to \mathfrak{gl}(\h)$ is the corresponding Lie algebra 
homomorphism. 
An IFAS $\nabla$ is said to be {\it reducible} if there exists a semidirect sum decomposition $\h =\da \ltimes \h'$ satisfying the following conditions:
\begin{enumerate}
  \item $\dim \da =1$,
  \item $g(\h') \da=0$,
  \item $\h'$ is invariant under the action of $g$, i.e. $g(\h)\h' \subset \h'$.
\end{enumerate}
\end{df}

If $g$ and $\h$ satisfies the above condition of reducible IFAS,  
then $g$ induces the homomorphism $g': \h \to \mathfrak{gl}(\h')$. Then we have $g'(X)Y = [X, Y]$ 
for $X \in \da$ and $Y \in \h'$. 
Furthermore the restriction of 
$g'$ to $\h'$ gives the homomorphism $g'|_{\h'}: \h' \to \mathfrak{gl}(\h')$, which corresponds to an IFAS on $\h'$. 

\begin{prop}\label{sum of IFAC2}
Suppose that a Lie algebra $\mathfrak{l}$ admits an $IFPS$, and a Lie algebra $\mathfrak{h}$ admits a reducible IFAS. Then $\mathfrak{l} \oplus \mathfrak{h}$ admits an IFAS.
\end{prop}
\pf 
We shall construct a Lie algebra homomorphism $h : \mathfrak{l} \oplus \mathfrak{h}$ $\to \mathfrak{gl}(\mathfrak{l} \oplus 
\mathfrak{h})$ corresponding to an IFAS as follows. Let $g : \h \to \mathfrak{gl}(\h)$  be a Lie algebra homomorphism corresponding to 
the given reducible IFAS on $\h$. We can decompose $\h$ into $\langle Z \rangle \ltimes \h'$ satisfying the condition of a 
reducible 
IFAS above. Let $g' : \h \to \mathfrak{gl}(\h')$ be the induced Lie algebra homomorphism.  
Note that the restriction $g'|_{\h'}$ satisfies the condition $(\ast)$ 
in Corollary \ref{IFAC and hom}. 

On the other hand $\l$ admits an IFPS from the assumption. Hence from Lemma \ref{lem3} and its remark, there exists a Lie algebra 
homomorphism $f : \l \oplus \mathfrak{a}_1 \to \mathfrak{gl}(\l \oplus \mathfrak{a}_1)$ satisfying $f(X)Y-f(Y)X = [X, Y]$ 
for $X$, $Y$ $\in \l \oplus \mathfrak{a}_1$, and 
$f(e) = I$ for the fixed basis $e$ of $\mathfrak{a}_1$.  
In the following we identify $Z$ with $e$. 
We define a representation $\phi: \l \oplus \da_1 \to \gl(\h')$ by $\phi(X)Y=0$, $\phi(e)Y = [Z, Y]$ for $X \in \l$ and $Y \in \h'$. 
Since $\phi$ is a derivation, we obtain the semidirect sum $(\l \oplus \da_1) \ltimes_\phi \mathfrak{h}'$. Then 
two Lie algebras $\mathfrak{l} \oplus (\langle Z \rangle \ltimes \mathfrak{h'})$ and $(\l \oplus \da_1) \ltimes_\phi 
\mathfrak{h}'$ are isomorphic, hence we may identify these two Lie algebras. 
By using $f$ and $g'$ we define a Lie algebra homomorphism 
$h:$ $\mathfrak{l} \oplus (\langle Z \rangle \ltimes \mathfrak{h'})$ $\rightarrow$ $\mathfrak{gl}((\l \oplus \da_1) 
\ltimes_\phi \mathfrak{h}')$ by
\begin{eqnarray*}
 h(X, 0) &=& \ 
\left(
\begin{array}{c|c}
f(X) & 0 \\
\hline
0 & 0
\end{array}
\right) \ ( X \in \l),\\
h(0, Z) &=& \ 
\left(
\begin{array}{c|c}
f(e)  & 0 \\
\hline 
0  & g'(Z)
\end{array}
\right), \\
h(0, Y) &=& \ 
\left(
\begin{array}{c|c}
0 & 0 \\
\hline
0 & g'(Y)
\end{array}
\right) \ (Y \in \h').
\end{eqnarray*}
Here the first row (column) in matrices stands for the $\l \oplus \da_1$-component, 
and the second row (column) stands for the $\h'$-component.
We can easily check that $h$ is a Lie algebra homomorphism. We show that $h$ satisfies the condition $(\ast)$ 
in Corollary \ref{IFAC and hom}: 
$h(X)Y-h(Y)X=[X,Y]$ for $X$, $Y$ $\in \mathfrak{l}\oplus \mathfrak{h}$. For two cases $X$, $Y$ $\in \l$, and $X$, $Y$ $\in \h'$, the 
condition $(\ast)$ is obviously satisfied. Hence 
we consider the remaining three cases:
\begin{enumerate}
  \item $X = Z$, $Y \in \l$,
  \item $X = Z$, $Y \in \h'$,
  \item $X \in \l$, $Y \in \h'$.
\end{enumerate} 
For $X \in \l \oplus \da_1$ and $Y \in \mathfrak{h}'$,  
we denote by
$
\begin{pmatrix}
X  \\
Y
\end{pmatrix}
$
the element of $(\l \oplus \da_1) \ltimes_\phi \mathfrak{h}'$. 
Then note that 
$
(0, Z) \in \mathfrak{l} \oplus (\langle Z \rangle \ltimes \mathfrak{h'})
$
is identified with 
$
\begin{pmatrix}
e  \\
0
\end{pmatrix}
$
$\in (\l \oplus \da_1) 
\ltimes_\phi \mathfrak{h}'$. 
In each case we have the following equality: \\[2mm]
(i)  
\begin{eqnarray*}
&{}& h(0, Z) \begin{pmatrix} Y \\ 0 \end{pmatrix}  - h(Y, 0) 
\begin{pmatrix} e \\ 0 \end{pmatrix} \\
&=& \begin{pmatrix} f(e)Y \\ 0 \end{pmatrix} - \begin{pmatrix} f(Y)e \\ 0 \end{pmatrix} = \begin{pmatrix} [e, Y] \\ 0 \end{pmatrix} 
= \begin{pmatrix} 0 \\ 0 \end{pmatrix} 
= [(0, Z), (Y, 0)].  
\end{eqnarray*}
(ii) 
\begin{eqnarray*}
&& h(0, Z) \begin{pmatrix} 0 \\ Y \end{pmatrix} - h(0, Y) 
\begin{pmatrix} e \\ 0 \end{pmatrix} \\
&=& \begin{pmatrix} I_{\l \oplus \langle Z \rangle} & 0 \\ 0 & g'(Z) \end{pmatrix}
\begin{pmatrix} 0 \\ Y \end{pmatrix} - \begin{pmatrix} 0 & 0 \\ 0 & g'(Y) \end{pmatrix} 
\begin{pmatrix} e \\ o \end{pmatrix} = 
\begin{pmatrix} 0 \\ g'(Z)Y \end{pmatrix} \\
&=& \begin{pmatrix} 0 \\ [Z, Y] \end{pmatrix}  = [(0, Z), (0, Y)]. 
\end{eqnarray*}
(iii) 
\begin{eqnarray*}
h(X, 0) \begin{pmatrix} 0 \\ Y \end{pmatrix} - h(0, Y) \begin{pmatrix} X \\ 0 \end{pmatrix} = 0 = 
[(X, 0), (0, Y)].
\end{eqnarray*}
\\[-2mm]

It follows that $\mathfrak{l}\oplus \mathfrak{h}$ admits an IFAS.
\hfill $\Box$\\
\renewcommand{\theenumi}{\alph{enumi}}
\renewcommand{\labelenumi}{\rm (\theenumi)}
\begin{prop}\label{decomposition}
Let $\l$ $=$ $\oplus \l_i$ be a decomposition of a Lie algebra $\l$ into its indecomposable ideals $\l_i$. If $\l$ satisfies the 
following condition {\rm(a)} or {\rm(b)}, then $\l$ admits an IFAS. \\[-4mm]
\begin{enumerate}
  \item Each $\l_i$ admits an IFAS.
  \item $\l_1$ admits an IFPS, $\l_2$ admits a reducible IFAS, and the remaining ideals $\l_i$ {\rm(}$i \neq 1, 2${\rm)} admit IFASs.
\end{enumerate}
\end{prop}
\pf 
In the case (a) its proof follows immediately from Proposition \ref{sum of IFAC1}. So let us consider the case (b). From 
Proposition \ref{sum of IFAC2},  $\l_1 \oplus \l_2$ admits an IFAS. Now we have ideals $\l_1 \oplus \l_2$ and $\l_i$ ($i \neq 1, 2$), 
each of which admits an IFAS. Since $\l$ is the direct sum of these ideals, $\l$ admits an IFAS. 
\hfill $\Box$

\newpage
\begin{prop}\label{decomposable case}
Any decomposable Lie algebra $\l$ of dimension $\leq 5$ admits an IFAS.
\end{prop}
\pf 
We decompose $\l$ into a direct sum of Lie algebras $\l = \oplus_{i=1}^m \l_i$ $(2 \leq m \leq 5)$, where $\l_i$ is an indecomposable 
ideal. Then from the  
proof of indecomposable case in Theorems \ref{IFPS} and \ref{IFAC} in $\S$5, each $\l_i$ admits an IFPS. 
If every ideal $\l_i$ admits an IFAS, then $\l$ admits an IFAS from Proposition \ref{decomposition}. Thus we suppose that 
there exists an ideal $\l_j$ which 
admits no IFAS. From the proof of indecomposable case in Theorem \ref{IFAC} in $\S$5, we have $\l_j =$ $\sll(2, \R)$ or 
$\mathfrak{o}(3)$.  
We may assume that $j = 1$.
Then there are only three possibilities of decomposition of $\l$; 
\[
\left\{
\begin{array}{l}
\l_1 \oplus A_{1,1}\\
\l_1 \oplus A_{1,1} \oplus A_{1,1} \\
\l_1 \oplus A_{2,1}. \\
\end{array}
\right.
\]
Recall that $A_{1,1}$ is the one dimensional abelian Lie algebra, and $A_{2,1}$ is the only 2-dimensional Lie algebra. 
From the remark of Lemma \ref{lem1}, $\l_1 \oplus A_{1,1}$ and  $\l_1 \oplus A_{1,1} \oplus A_{1,1}$ admit an IFAS. 
From the proof of indecomposable case in Theorem \ref{IFAC} in $\S$5, 
$A_{2,1}$ admits an IFAS $\nabla$. Considering the construction of 
Lie algebra homomorphism in Proposition \ref{solv and IFAC}, it is easy to see that
this IFAS $\nabla$ is reducible with respect to the
decomposition of $A_{2,1} = \langle X_1 \rangle \ltimes \langle X_2 \rangle$. Hence 
$\l_1 \oplus A_{2,1}$ admits an IFAS by Proposition \ref
{sum of IFAC2}. 
\hfill $\Box$\\

By the above discussions Theorems \ref{IFPS} and \ref{IFAC} have been completely proved.\\[-1mm]

\remark
By virtue of Proposition \ref{sum of IFAC2} we can prove Proposition \ref{decomposable case} without constructing  
representations of $\sll(2, \R) \oplus A_{2,1}$ and $\mathfrak{o}(3) \oplus A_{2,1}$ corresponding to IFASs. 
There is another purpose for proving Proposition \ref{sum of IFAC2}. It is a construction of Lie groups of higher dimension 
admitting an IFAS.  
For example in $\S 5$ we showed that any non-perfect indecomposable Lie algebra $\h$ of dimension $\leq 5$ admits an IFAS. 
We can easily see that it is reducible. 
Hence for any Lie algebra $\l$ admitting an IFPS, $\l \oplus \h$ admits an IFAS by Proposition \ref{sum of IFAC2}.  

There is another example of a Lie algebra admitting a reducible IFAS. It is the set of upper triangular matrices $\mathfrak{t}(n, \R)$. This 
Lie algebra 
$\mathfrak{t}(n, \R)$ decomposes into a semidirect sum $\h \ltimes \kk$, where $\h$ is the set of diagonal matrices and $\kk$ is the 
derived subalgebra. We graduate $\kk$ naturally, i.e. we set $\kk_i$ $=$ 
$\langle e_{1, i+1}, e_{2, i+2}, \ldots,$ $e_{n-i,  n} \rangle$, 
where $e_{i, j}$ is the matrix having $1$ in the $(i, j)$ position and $0$ elsewhere.  
The decomposition $\h \ltimes \oplus_{i} \kk_i$ satisfies the condition of Proposition \ref{solv and IFAC}. Hence 
the Lie algebra  
$\mathfrak{t}(n, \R)$ admits a representation $g:$ $\mathfrak{t}(n, \R) \to \gl(\mathfrak{t}(n, \R))$ corresponding to an IFAS.   
We denote by $\mathfrak{a}$ the one-dimensional subspace spanned by $I_n$ $\in$ $\h$, and denote by $\mathfrak{b}$ its arbitrary 
complementary subspace in $\h$. 
Hence we have $\h = \mathfrak{a} \oplus \mathfrak{b}$.   
From the definition of $g$ in Proposition \ref{solv and IFAC} and the fact that $\h$ is abelian, we can see the following easily: 
the IFAS corresponding to $g$ and the decomposition $\mathfrak{t}(n, \R) = \mathfrak{a} \ltimes (\mathfrak{b} \ltimes 
\kk)$ satisfies the condition of reducible IFAS.   
Hence for any Lie algebra $\l$ admitting an IFPS, the direct sum $\l \oplus \mathfrak{t}(n, \R)$ admits an IFAS. 

\section*{Acknowledgement}
The author expresses his sincere thanks to Prof. Y. Agaoka who led him to study the subject and read through the manuscript in 
detail.

\begin{small}
  \parindent0pt\parskip\baselineskip

  Hironao Kato \\
  Department of Mathematics, Graduate School of Science, Hiroshima University, Higashi-Hiroshima 739-8526, Japan \\ 
  
  \textit{E-mail}: katoh-in-math@hiroshima-u.ac.jp
\end{small}
\end{document}